\documentstyle[amscd]{amsart}
\newtheorem{thm}{Theorem}[section]
\newtheorem{lm}[thm]{Lemma}
\newtheorem{cor}[thm]{Corollary}
\newtheorem{pro}[thm]{Proposition}

\newtheorem{th}{Theorem}

\theoremstyle{definition}
\newtheorem{df}[thm]{Definition}

\theoremstyle{remark}
\numberwithin{equation}{section}

\def \BR {\Bbb{R}}
\def \BZ {\Bbb{Z}}
\def \BQ {\Bbb{Q}}
\def \CA {\cal A}
\def \CB {\mathcal B}

\def \CG {\cal G}

\def \CM {\cal M}

\def \CR {\cal R}

\def \hM {\hat{{\cal M}}}

\def \la {\langle}
\def \ra {\rangle}
\def \a {\alpha}
\def \b {\beta}

\def \d {\delta}

\def \Lam {\Lambda}

\def \o {\omega}

\def \Si {\Sigma}

\def \bd {\partial}

\def \x {\times}
\def \ox {\otimes}

\begin{document}

\baselineskip.525cm

\title[Cap-product structure]
{Cap-product Structures \\ on the Fintushel-Stern spectral sequence}

\author[Weiping Li]{Weiping Li}
\address{Department of Mathematics, Oklahoma State University \newline
\hspace*{0.175in}Stillwater, Oklahoma 74078-0613}
\email{wli@@math.okstate.edu}
\thanks{Partially supported by NSF Grant DMS9626166. }
\date{October 7, 1997.}

\begin{abstract}
We show that there is a well-defined cap-product structure on the
Fintushel-Stern spectral sequence. Hence we obtain the induced
cap-product structure on the ${\BZ}_8$-graded instanton Floer homology.
The cap-product structure provides an essentially new property of the instanton
Floer homology, from a topological point of view, which 
multiplies a finite dimensional cohomology class by an infinite dimensional 
homology class (Floer cycles) to get another infinite dimensional homology
class.
\end{abstract}

\maketitle

\section{Introduction}

In \cite{at} p298, Atiyah stated that {\em In the product we then have four
types of cycle
\[\text{finite} \x \text{finite}, \hspace{0.6in}
\text{cofinite} \x \text{cofinite},\]
\[\text{cofinite} \x \text{finite}, \hspace{0.6in}
\text{finite} \x \text{cofinite}.\]
The first two give ordinary homology and cohomology respectively. The
other two are quite different from these and give in an obvious sense
``middle-dimensional'' homology, one ``positive'' and one ``negative.''}
It is the purpose of this paper to elaborate on this comment and to
show this point rigorously.
A. Floer \cite{fl} introduced a novel homology theory for oriented
closed 3-manifolds with the homology groups of $S^3$. Subsequently, 
Donaldson \cite{at} discovered that the instanton Floer homology has an
intimate relationship to his gauge theoretic 4-manifold invariants \cite{dk}.
More precisely, if $X= X_1 \cup_Y X_2$, where $X_i (i=1, 2)$ is a
4-manifold with homology 3-sphere $Y$ boundary, then
the relative Donaldson invariants $\Phi (X_1): H_2(X_1; {\BZ}) \to 
HF_*(Y)$ and $\Phi (X_2): H_2(X_2; {\BZ}) \to HF_*(\overline{Y})$ lie in
the Floer homology cycles ($\overline{Y}$ has opposite orientation of $Y$).
The natural pairing in the Floer theory provides
$\Phi (X) = \la \Phi (X_1), \Phi (X_2) \ra$ (see \cite{at, br, bd}).
The pairing of the Floer cycles defines 4-dimensional Donaldson invariant.

In a formal way, the Floer homology can be viewed as the 
``middle-dimensional'' homology of the infinite dimensional space
${\CB}_Y$ of gauge equivalence classes of $SU(2)$-connections on $Y \x SU(2)$.
Let $Z(X_i) (i=1, 2)$ be the set of boundary values of instantons over 
$X_i$. The $Z(X_i)$ represents a cycle of homology class $\Phi (X_i)$
which is independent of metrics.
The pairing $\la \Phi (X_1), \Phi (X_2) \ra$ can be identified with the
intersection $Z(X_1) \cap Z(X_2)$ in ${\CB}_Y$.
The formalization alone might suggest that there should be a product from
$H^*({\CB}_Y) \ox HF_*(Y)$ to $HF_*(Y)$ (see \cite{bd} \S 4.2), as multiplying
a finite dimensional cohomology class
by an infinite dimensional class to get another infinite dimensional class
(see \cite{at} p.298). Such a cap-product structure on the
instanton Floer homology shows that the Floer homology groups are, from a 
topological point of view, something essentially new \cite{at} p298.
Dually, this is equivalent to a cup-product structure on the instanton
Floer cohomology $H^*({\CB}_Y) \ox HF^*(Y) \to HF^*(Y)$ (see \cite{at, bd}).

In this paper, we show that there is a well-defined
cap-product structure of $H^*({\CB}_Y, {\BQ})$ on the ${\BZ}_8$-graded
instanton Floer homology. For a single 3-dimensional cohomology class in
$H^*({\CB}_Y, {\BQ})$, this cap-product 
is known \cite{br, bd}. In order to obtain the 
cap-product structure, we have to overcome the main difficulty arising from
the {\it non-compactness} of the moduli space
${\CM}_{Y \x {\BR}}(a, b)$ with flat connections $a$ and $b$ over $Y$.
For a cohomology class $\o \in H^*({\CB}_Y, {\BQ})$ with degree $\geq 8$, 
the pairing $\la \o, {\CM}_{Y \x {\BR}}(a, b)\ra$ in general is not defined,
since ${\CM}_{Y \x {\BR}}(a, b)$ has three possible non-compactness cases
from (1) the ${\BR}$-action, (2) the concentrated instanton bubbling on
$Y \x {\BR}$ and (3) the chain anti-self-dual connection splitting the moduli
space into codimension one (or $\geq 1$) pieces. The non-compactness of the 
${\BR}$-action is not serious. We can work on the space
${\CM}_{Y \x {\BR}}(a, b)/{\BR} = {\hM}_{Y \x {\BR}}(a, b)$.
The non-compactness of the bubbling type (2) is problematic.
We use the ${\BZ}$-graded instanton Floer homology defined by Fintushel-Stern
\cite{fs} to fix Chern-Simons values in one of the energy band
$(r, r+1)$, where $r$ is a regular value of the Chern-Simons functional.
Any bubbling required energy at least 1 forces the corresponding
moduli space to have different lifting for the generators $a$ and $b$.
So the moduli space with fixed lifting $a^{(r)}$ and $b^{(r)}$ contains 
no bubbling limit point. 
The {\it chain anti-self-dual connections} also provide the codimension one
components in the compactification. The evaluation of cohomology class
on the compactification does not make sense. 
This is why we restrict ourself to the rational cohomology
group $H^*({\CB}_Y, {\BQ})$. The rational cohomology group $H^*({\CB}_Y, {\BQ})$
is generated by the $\mu$-map construction. Thus we can represent each
cohomology class by its Poincar\'{e} dual as a divisor from the homology
cycle of the 3-manifold $Y$ itself. Thus adapting the 4-dimensional technique
(dimension counting replaced by spectral flow counting), we can
show that the set $P.D(\o) \cap {\hM}_{Y \x {\BR}}(a, b)$ is compact
by ruling out those three non-compact parts.
In turn, there is a spectral sequence which starts
from the ${\BZ}$-graded Floer homology and converges to the
${\BZ}_8$-graded Floer homology. Thus the cap-product structure on the
${\BZ}$-graded Floer homology induces the cap-product structure on the
spectral sequence, hence on the ${\BZ}_8$-graded Floer homology.
Our main theorem is following.

\noindent{\bf Theorem} [Theorem~\ref{main}]
{\em There is a well-defined cap-product action of $H^*({\CB}_Y; {\BQ})$ on the
Fintushel-Stern spectral sequence $(E^k_{n,j}(Y), d^k)$.}

In particular, Theorem~\ref{main} shows that there is a well-defined cap-product
(cup-product) structure on the instanton Floer homology (cohomology). There
are some cross-product structure on the $SO(3)$ Floer homology which proved
by Braam and Donaldson for rational coefficients and by the author
for integral coefficients (see \cite{li2} for more references).
Note that the space ${\CB}_Y$ is not globally a product but only
infinitesimally (because its tangent bundle decomposes). {\em
Thus its middle-dimensional homology cannot be reduced to ordinary homology and
cohomology by a factorization. Thus the Floer homology groups
are, from a topological point of view, something essentially new} (quoted from
\cite{at} page 298).

The paper is organized as follows. In \S 2, we give a review on
the instanton Floer theory of integral homology 3-spheres.
We start the rational cohomology group and ring of ${\CB}_Y$ in \S 3.1, and
show that there is a well-defined cap-product action on the
${\BZ}$-graded instanton Floer homology in \S 3.2. Using the 
commutative property of the cap-product action with higher differentials, 
we prove our main theorem in \S 4.

\section{The instanton Floer homology of integral homology 3-spheres}
 
In this section, we give a description of the gauge theory on
3-manifolds and review the definition of the instanton Floer homology (see
\cite{br, fs, fl}).
\vspace{0.1 in}

Let $Y$ be an oriented closed 3-dimensional
smooth manifold with $H_1(Y, {\BZ}) = 0$, and let $P \rightarrow
Y$ be a smooth principal $SU(2)$-bundle (always trivial). 
Fix a trivialization $Y\times SU(2)$ of $P$ and let $\theta $ 
be the associated trivial connection.
Denote the Sobolev $L_k^p$ space of connections on $P$ by ${\CA}_Y$. 
This space has a
natural affine structure with underlying vector space ${\Omega
}^1(Y, adP)$ where $adP$ is the adjoint bundle. ${\CA}_Y$ is acted
upon by the gauge group ${\CG} = L_{k+1}^p({\Omega }^0(Y,adP))$ 
of bundle automorphisms of $P$, and the orbit space ${\CB}_Y = {\CA}_Y
/{\CG}$ is well-defined when $k+1 > 3/p$. The universal cover
$\tilde{{\CB}}_Y$ of ${\CB}_Y$ defined as $\tilde{{\CB}}_Y =
{\CA}_Y/{\CG}_0$, where ${\CG}_0 \subset {\CG}$ is the group
of degree zero gauge transformations. The irreducible connections form
an open and dense subspace ${\CB}^{\ast }_Y$ of ${\CB}_Y$ which
is a Banach manifold with
\[ T_a {\CB}^{\ast }_Y \equiv  \{ b \in L_k^p ({\Omega
}^1(Y, adP)) | \ \ d_a^{\ast } b = 0 \} \]
where $d_a^{\ast }$ is the $L^2$-adjoint of $d_a$ (covariant derivative 
on sections of $adP$) with respect to some metric on $Y$.
\vspace{0.1 in}

The Chern-Simons functional $cs: {\CA}_Y \rightarrow {\BR}$ is defined
as \[ cs(a) = \frac{1}{8\pi^2} \int_Y tr(a\wedge da + \frac{2}{3}a\wedge
a\wedge a).\] and  satisfies
$ cs(g\cdot a) = cs(a) + \deg(g)$ 
for gauge transformations $g:Y \rightarrow SU(2)$.
Thus $cs$ is well-defined on $\tilde{{\CB}}(P) = {\CA}(P)/\{g \in
{\CG}: \deg(g)=0 \}$  and it descends to a function
$cs : {\CB}(P) \rightarrow {\BR}/{\BZ}$
which plays the role of a Morse function in defining
Floer homology. The differential of $cs$ is
\[dcs(a)(\a ) = \int_Y tr(F_a \wedge \a),\]
hence its critical set consists of the flat connections
${\CR}({\CB}(P)) =\{ a \in {\CB}(P) |\ \  F_a = 0 \}$ (I.e., $F_a$ is
the curvature 2-form on $Y$). It is well-known that the elements of
${\CR}({\CB}(P))$
are in 1-1 correspondence with those of 
\[{\CR}(Y)= Hom(\pi_1(Y),SU(2))/{ad
SU(2)}, \]
the $SU(2)$-representations of $\pi_1(Y)$ modulo conjugacy.
Let ${\CR}^*(Y) = {\CR}(Y) \setminus \{\theta \}$ be the space of 
irreducible $SU(2)$ representations.

Define the weighted Sobolev space $L_{k,\delta
}^p$ on sections $\xi $
of a bundle over $Y\times {\BR}$ to be the space of $\xi$ for
which $e_{\delta }\cdot \xi $ is in $L_k^p$,
where $e_{\delta }(y,t) = e^{\delta |t|}$ for
$|t| \geq 1$. For any $\delta \geq 0$ and any $SO(3)$ anti-self-dual
connection $A$ on the bundle over $Y \times {\BR}$,
the anti-self-duality operator
\[d_A^{\ast } \oplus d_A^+ : L_{k+1,\delta }^p({\Omega }^1(Y\times {\BR},
adP)) \rightarrow
L_{k,\delta }^p(({\Omega }^0 \oplus {\Omega }^2_+)(Y\times {\BR},adP))\]
is Fredholm. We call $A$  {\em regular} if $d_A^{\ast } \oplus d_A^+ $
is surjective and we call the moduli
space ${\CM}_{Y \times {\BR}}$
of perturbed anti-self-dual connections with finite energy regular
if it contains only orbits of regular $A$'s.
The spectral flow
is $SF(a_{\a }, a_{\b }) = Index(d_A^{\ast } \oplus d_A^+ )(a_{\a }, a_{\b })
\pmod {2 p_1(A)}$, where $p_1(A)$ is the Pontryagin form. 
The $Index(d_A^{\ast } \oplus d_A^+ )(a_{\a }, a_{\b })$ is given by
\begin{equation} \label{Atp}
Index(d_A^{\ast } \oplus d_A^+ )(a_{\a }, a_{\b })= -2 \int_{Y \times {\BR}}
p_1(A) - \frac{h_{\beta } + \rho_{\beta }(0)}{2} + \frac{ - h_{\a } +
\rho_{\a }(0)}{2}, 
\end{equation}
where 
$h_{\a }$ is the sum of the dimensions of 
$H^i(Y, ad a_{\a }) (i =0, 1)$,
and $\rho_{\a }(0)$ is the $\rho $-invariant of the signature operator
$* d_{a_{\a }} - d_{a_{\a }}*$ over $Y$ restricted to even forms (see
\cite{fs, fl}). We use the trivial connection and a
fixed tangent vector to the trivial connection to fix $\mu (\theta ) = 0$.
So
 \[ \mu (a_{\a }) \equiv Index (d_A^{\ast } \oplus
 d_A^+)(a_{\a } ,\theta ) \ \ \  (\mbox{mod}\ \  8) . \]
Floer's chain group $ C_j(Y)$ is defined to be the free module
 generated by irreducible flat 
 connections $a_{\a } \in f^{-1}(0)$ 
 with $\mu (a_{\a }) \equiv j \pmod 8$.

The Chern-Simons and spectral flow functionals on ${\CA}_Y$
descends to functionals $cs$ and $SF: \tilde{{\CB}}_Y \to {\BR}$.
Let $cs ({\CR}_Y^*)$ be the image values of Chern-Simons functional
in ${\BR}$. So $cs({\CR}_Y^*) \pmod 1$ is a finite set and independent of
the choice of $\theta $. The set ${\BR}_Y = {\BR}
\setminus cs({\CR}_Y^*)$
consists of all the regular values of the Chern-Simons functional.
Given $a \in {\CR}^*_Y \subset {\CB}_Y$, let $a^{(r)} \in
\tilde{{\CR}}_Y \subset \tilde{{\CB}}_Y$ be the unique lift of a
such that $cs(a^{(r)}) \in (r, r+1)$ for $r \in {\BR}_Y$.
Note that $SF(a^{(r)}) = SF(a^{(r)}, \theta)$ is an integer.
 
For $a, b \in {\CB}_Y$, choose any smooth representatives $a, b \in
{\CA}_Y$ and a connection $A$ on $Y \times {\BR}$
whose $t$-component vanishes and which equals $a$ for large negative $t$
and equals $b$ for large positive $t$.
Set ${\CA}_{\delta }(a, b) = A + L^4_{1, \delta}(\Omega^1
(Y\times {\BR},  adSU(2))$. It is acted smoothly upon by
the gauge group
\[{\CG}_{\delta } = \{g\in L^4_{2,\mbox{loc}}(Y\times {\BR},
adSU(2)) | \ \ \exists T > 0, \]
\[ \xi \in L^4_{2, \delta
}(\Omega^0 (Y\times {\BR}, adSU(2))) \ \ \ \mbox{such that} \ \ \
g= \exp (\xi) \ \
\mbox{for} \ \ \ |t| \geq T \}. \]
Then the quotient space ${\CB}_{\delta }(a, b) = {\CA}_{\delta
}(a,b)/{\CG}_{\delta }$ is a smooth Banach manifold.
Note that the space ${\CB}_{\delta }(a, b)$ depends on the homotopy
class of path $A$ of connections between $a$ and $b$.
Each $A \in {\CA}_Y(a, b)$
is (temporal) gauge equivalent to a connection whose
component in the ${\BR}$-direction vanishes.
Floer's crucial observation is that trajectories of the vector field
$f$, i.e., the flow lines of \[\frac{\bd a}{\bd t} + f(a(t)) =
0, \ \ \ \mbox{or} \ \ \ \
\frac{\bd a}{\bd t} =  \star F(a(t)),\]  can be
identified with instantons $A$ on $Y \times {\BR}$
and $A|_{Y\times \{t \}} = a(t)$. 
Let ${\CM}(a_{\a }, a_{\beta })$ be the 
space of (perturbed) anti-self-dual connections 
$A$ with asymptotics  
flat 
connections $a_{\a }$ and $a_{\beta }$.
The space ${\CM}(a_{\a }, a_{\beta })$ consists of 
infinitely many smooth manifold components,
the dimensions of various components differ by a
multiple of $8$ and all equal to 
$\mu (a_{\a }) - \mu (a_{\beta }) \pmod 8$.
Let ${\CM}^1(a_{\a }, a_{\beta })$
be the union of all the 1-dimensional components 
in ${\CM}(a_{\a }, a_{\beta
})$. Then ${\CM}^1(a_{\a }, a_{\beta })$ has a canonical orientation
together with a proper, free ${\BR}$-action induced by the translations
on $Y \times {\BR}$. 
It follows that ${\CM}^1(a_{\a }, a_{\beta })/{\BR} = \hat{{\CM}}(
a_{\a }, a_{\beta })$ is a compact oriented
0-manifold
(a finite set
of signed
points) by the Proposition 1c.2 in \cite{fl}.
The differential $\bd_Y: C_j(Y) \rightarrow C_{j-1}(Y)$ of the Floer chain
complex is defined by
\[ \bd_Y a_{\a } = \sum_{a_{\beta} \in C_{j-1}(Y)} \# \hat{{\CM}}(a_{\a }
,a_{\beta} )a_{\beta} , \]
where $ \# \hat{{\CM}}(a_{\a }, a_{\beta })$ is the 
algebraic number of points, and the sign is  
given by the spectral flow.
Floer has shown that ${\partial }^2 = 0$, hence 
the homology of this
complex $\{C_j(Y), \bd_Y\}_{j \in {\BZ}_8}$
is the (instanton) Floer homology $HF_*(Y) (* \in {\BZ}_8)$.
It is shown that $HF_*(Y)$ is independent
of the choice of metric $\sigma $ on $Y$ and of regular   
perturbations by the Theorem 2 in \cite{fl}.

For lifts $\tilde{a}$ and $\tilde{b}$ in $\tilde{\CB}_Y$, denote
$\tilde{\CB}_{\d}(\tilde{a}, \tilde{b})$ for lifts of connections
in ${\CB}_{\d}(a, b)$.
If $A \in {\CB}_{\d}(a, b)$ for $a, b \in {\CR}^*_Y$
and $\tilde{A}$ is any lift to $\tilde{{\CB}}_{\delta
}(\tilde{a}, \tilde{b})$, then ind$D_{\tilde{A}} = SF(\tilde{a}) -
SF(\tilde{b})$.
For $\mu (a) \equiv j  \pmod 8 , b = \theta $ and $r \in {\BR}_Y$, 
the spectral flow $SF(a^{(r)})$ is independent of the choice of trivial
connection $\theta $ used to define Chern-Simons functional on 
${\CA}_{\delta }(a, \theta)$.

\begin{df} Define the chain groups $C_n^{(r)}(Y) = {\BZ}
 \{a \in {\CR}^*_Y \ \ 
| \ \ SF(a^{(r)}) = n \in \BZ \}$
and $\bd^{(r)} : C_n^{(r)}(Y) \to C_{n-1}^{(r)}(Y)$
 by 
 \[ \bd^{(r)}a = \sum_{b\in C_{n-1}^{(r)}(Y)} \# \hat{{\CM}}^1_{Y
 \times {\BR}}(a, b) \cdot b . \]
\end{df}

Note that $\bd^{(r)}a$ 
defined in \cite{fs} can be identified with
$\bd^{(r)}a = \sum_{b\in C_{n-1}^{(r)}(Y)} \# \hat{{\CM}}^1_{Y
 \times {\BR}}(a^{(r)}, b^{(r)}) \cdot b$. For any $A \in \hat{{\CM}}^1_{Y
 \times {\BR}}(a, b)$, there is a unique lift $\tilde{A} \in
\hat{{\CM}}^1_{Y  \times {\BR}}(a^{(r)}, \tilde{b})$ (similarly for
$\tilde{A} \in \hat{{\CM}}^1_{Y  \times {\BR}}(\tilde{a}, b^{(r)})$. Since
$Ind_{\tilde{A}} = SF(a^{(r)}) - SF (\tilde{b}) = 1$, we have
$SF (\tilde{b}) = n-1$ and $\tilde{b} = b^{(r)}$ the preferred lift.
Hence $\tilde{A} \in \hat{{\CM}}^1_{Y  \times {\BR}}(a^{(r)},b^{(r)})$.
The orientation of $A$ and $\tilde{A}$ are same, so the
counting $\# \hat{{\CM}}^1_{Y \times {\BR}}(a, b)$ can be also
expressed as $\# \hat{{\CM}}^1_{Y \times {\BR}}(a^{(r)}, b^{(r)})$.

The key observation is that the Chern-Simons functional is non-decreasing along 
the gradient trajectories in $\tilde{{\CB}}_Y$.
By the similar argument of Floer \cite{fl}, 
Fintushel and Stern \cite{fs} obtained: 
\begin{enumerate}
\item $\bd_{n-1}^{(r)} \circ \bd_n^{(r)} = 0$ for $n \in {\BZ}$.
The homology of $\{ C_n^{(r)}(Y) , \bd_n^{(r)}\}$, denoted by
$H_*(C_n^{(r)}(Y) , \bd_n^{(r)}) = I_*^{(r)}(Y), * \in {\BZ}$, is 
the ${\BZ}$-graded Floer homology which is independent of metrics and
perturbations.
\item If $[r, s] \subset {\BR}_Y$, then $I_*^{(r)}(Y) = I_*^{(s)}(Y)$ and
$I_*^{(r+1)}(Y) = I_{* +8}^{(r)}(Y)$.
\item The ${\BZ}$-graded Floer homology $\{I_n^{(r)}(Y), n \in {\BZ}\}$
determines the ${\BZ}_8$-graded Floer homology groups by filtering the Floer
chain complex.
\end{enumerate}

For $r \in {\BR}_Y, j \in {\BZ}_8, n \in {\BZ}$ and $n \equiv j \pmod 8 $,
define the free abelian groups
\begin{equation} \label{fil}
F_n^{(r)} C_j(Y) = \sum_{k \geq 0} C_{n + 8k}^{(r)}(Y). \end{equation}
Then there is  a finite length decreasing filtration of $C_j(Y) (j \in
{\BZ}_8)$:
\[\cdots \subset F_{p+8}^{(r)}C_j(Y) \subset F_p^{(r)}
C_j(Y) \subset F_{p-8}^{(r)}C_j(Y)
\subset  \cdots \subset C_j(Y) , \]
\[ C_j(Y) = \bigcup_{k \in {\BZ}}F_{j+8k}^{(r)} C_j(Y) .\]
Since the perturbed Chern-Simons functional is non-decreasing along
the gradient flows, it follows that Floer's
boundary map $\bd_Y: F_n^{(r)}C_j(Y) \to F_{n-1}^{(r)}C_{j-1}(Y)$ 
preserves the
filtration
\[ \begin{array}{ccccc}
 & \downarrow & \downarrow &  & \downarrow \\
\cdots & \subset F_{n+8}^{(r)}C_j(Y) & \subset F_n^{(r)}C_j(Y) & \subset
\cdots & \subset C_j(Y) \\
 & \downarrow \bd_Y & \downarrow \bd_Y &  & \downarrow \bd_Y \\
\cdots & \subset F_{n+7}^{(r)}C_{j-1}(Y) & \subset F_{n-1}^{(r)}C_{j-1}(Y) 
& \subset \cdots & \subset C_{j-1}(Y) \\
& \downarrow & \downarrow &  & \downarrow  
\end{array} . \]
So the ${\BZ}_8$-graded Floer chain complex has a decreasing
  bounded filtration $(F_n^{(r)}C_*(Y), * \in {\BZ}_8 )$. 
The homology of the vertical chain subcomplex in the above filtration gives
$F_n^{(r)}I_j(Y), $
\[ \begin{array}{ccccc} 
 & \downarrow & \downarrow &  & \downarrow \\
 \cdots & \subset F_{n+8}^{(r)}I_j(Y) & \subset F_n^{(r)}I_j(Y) & \subset
 \cdots & \subset HF_j(Y) \\ 
  & \downarrow  & \downarrow  &  & \downarrow  \\
  \cdots & \subset F_{n+7}^{(r)}I_{j-1}(Y) & \subset
  F_{n-1}^{(r)}I_{j-1}(Y) 
  & \subset \cdots & \subset HF_{j-1}(Y) \\
  & \downarrow & \downarrow &  & \downarrow  
  \end{array} . \]
This is a bounded filtration for the ${\BZ}_8$ graded Floer homology groups
$\{HF_j(Y) \}_{j \in {\BZ}_8}$:
\[\cdots \subset F_{n+8}^{(r)}HF_j(Y) \subset F_n^{(r)}HF_j(Y) 
\subset F_{n-8}^{(r)}HF_j(Y)
  \subset \cdots \subset HF_j(Y) , \]
where $F_n^{(r)}HF_j(Y) = Im[ F_n^{(r)}I_j(Y) \to HF_j(Y)]$.
By the definition of the filtration (\ref{fil}), 
\[F_n^{(r)}C_j(Y)/F_{n+8}^{(r)}C_j(Y) = C_n^{(r)}(Y)  , \]
and the induced chain map is precisely $\partial_n^{(r)}$.

\begin{thm} \label{E1} [Theorem 5.1 in \cite{fs}] 
(i) There is a spectral 
sequence $(E_{n,j}^k, d^k)$ (periodic in $j$ with period $8$) with 
\[ E_{n,j}^1(Y) \cong I_n^{(r)}(Y), \ \ \ \ 
d^k : E^k_{n,j}(Y) \to E^k_{n+8k-1,j-1}(Y), \]
\[E^{\infty}_{n,j}(Y) =
F_n^{(r)}I_j(Y)/F_{n+8}^{(r)}I_j(Y),
\ \ \ \ (j \in {\BZ}_8, n \in {\BZ}, n \equiv j \pmod 8 ). \] 
Furthermore, the groups $E^k_{n,j}(Y)$ are topological invariants. 

(ii) Any one period of the spectral sequence $(E_{n,j}^k(Y), d^k)$
converges to $E^{\infty}_{*,j}(Y) (j \in {\BZ}_8)$ which is isomorphic to
the bigraded module associated to the filtration $F^{(r)}_n$ of
$I_j(Y)\, (j \in {\BZ}_8)$:
\[ HF_j(Y) \cong \bigoplus_{k\in {\BZ}}E^{\infty}_{n +8k,j}(Y).\]
\end{thm} \qed

Note that our formulation of the spectral sequence is slightly different
from the one in \cite{fs}. We span the spectral sequence formulated in \cite{fs}
periodically in $j$-direction. So all higher differentials $d^k (k \geq 1)$
maps into the same direction with degree $(8k-1, -1)$.

\begin{cor} \label{ident}
For $j \in {\BZ}_8$, $\bigoplus_{k\in {\BZ}} I_{j+8k}^{(r)}(Y) = HF_j(Y)$ if
and only if $d^k = 0 (k \geq 1)$ in the spectral sequence.
\end{cor} 
Proof: Since the following is true for one period of the $E^{\infty}_{*, *}$:
\[HF_j(Y) = \bigoplus_{k\in {\BZ}}E^{\infty}_{n +8k, j}(Y), \]
and $E^1_{n,j}(Y)$ is given by the ${\BZ}$-graded Floer homology, so
the result follows from the fact that
$(E^k_{*,*}(Y), d^k)$ collapses at the $E^1$-term. \qed

\section{A cap-product structure on $I^{(r)}_*(Y, {\BZ})$}

In this section, we show that the ${\BZ}$-graded Floer homology 
$I^{(r)}_*(Y, {\BZ})$ admits a $H^*({\CB}_Y, {\BQ})$-cap-product structure: 
$H^*({\CB}_Y, {\BQ}) \ox I^{(r)}_*(Y, {\BZ}) \to I^{(r)}_*(Y, {\BZ})$. 
The main
task is to overcome the non-compactness when
one evaluates cohomology classes on the moduli spaces. This
is the essential point that we use the ${\BZ}$-graded Floer homology to
filter out the ideal instantons over $Y \x {\BR}$, and use
the rational cohomology
classes to rule out the chain connections (codimension one components).

\subsection{The cohomology ring $H^*({\CB}_Y, {\BQ})$}

In this subsection we describe that the cohomology ring (group) of
${\CB}_Y$ with rational coefficients is generated by two special classes.
Both classes can be obtained by the $\mu$-map procedure with divisors
realized by time-translation invariant divisors.

\begin{pro} \label{1e}
For an integral homology 3-sphere $Y$ and $a, b \in {\CR}^*(Y)$, the
rational cohomology ring $H^*(\tilde{\CB}_{Y \x {\BR}}(a, b); {\BQ})$
is the tensor product of a polynomial algebra on the four-dimensional
generator $\nu$ and an exterior algebra on the one-dimensional generator
$\mu (Y)$:
\[H^*(\tilde{\CB}_{Y \x {\BR}}(a, b); {\BQ}) =
\text{Sym}^*(H_0(Y \x {\BR})) \ox \Lambda^*(H_3(Y \x {\BR})).\]
\end{pro}
Proof: Since $H_*(Y, {\BZ}) = H_*(S^3, {\BZ})$, we have 
$H_1(Y\x {\BR}) = H_2(Y\x {\BR}) = H_4(Y\x {\BR}) = 0$;
$H_0(Y\x {\BR}) = H_3(Y\x {\BR}) = {\BZ}$. 
The rational cohomology of $\tilde{\CB}_{Y \x {\BR}}(a, b)$ is the product of
a polynomial algebra on even-dimensional generators and an
exterior algebra on odd-dimensional ones. The generators arise from
the $\mu_{- \frac{1}{4}p_1}$ construction for $\tilde{\CB}_{Y \x {\BR}}(a, b)$
(see \cite{dk} Chapter 5 and \cite{bd}). Here $p_1$ is the first 
Pontryagin class of the base-point fibration 
$\tilde{\CB}^*_{Y \x {\BR}}(a, b) \to {\CB}^*_{Y \x {\BR}}(a, b)$.
So the result follows. 
\qed

These two cohomology classes $\nu = \mu ((y_0, t_0))$ 
and $\mu(Y)$ are obtained through the
general universal $\mu$-map construction. We give
their geometric construction to pass down to the cohomology classes in
${\CB}_Y$.

For the purpose of finding the cap-product structures on the instanton Floer
homology, we can fix our asymptotic values as irreducible flat
connections corresponding elements in ${\CR}^*(Y)$. For 
$a$ and $b$ in ${\CR}^*(Y)$,
${\CB}^*_{Y \x {\BR}}(a,b)$
is the equivalence classes of irreducible
connections over $Y \x {\BR}$ with asymptotic values $a$ and $b$. Pick
a point in $Y$, call $y = \{y_0\}$, 
${\CG}_0$
is gauge transformation group of $SU(2)$-vector bundle which acts as
identity between the fibers over $\{y, + \infty\}$ and $\{y, -
\infty\}$ in $Y \x {\BR}$. This is a normal subgroup of
${\CG}_{\delta }$, thus we can form $\tilde{\CB}_{Y \x {\BR}}^*(a, b)
= {\CA}_{\delta }^* /{\CG}_0$. Let us choose a path
$\{y\} \x {\BR}$ in $Y \x {\BR}$
representing $[\BR]$. For each connection $A \in {\CA}_{\delta
}^*(a,b)$, we first fix the trivializations at $P_{\{y,\pm \infty\}}$
and identify them. Choose $p \in P_{\{y, - \infty\}}$ lies over $\{y,
- \infty \}$, then there is a unique path that is always horizontal
with respect to the connection $A$. So the parallel transport along
$\{y\} \x {\BR}$ gives an automorphism of the fibers between
$P_{\{y,- \infty\}}$ and $P_{\{y,+ \infty \}}$ and after a choice of
$p$ this determines an element of $SO(3)$ the honolomy of the path.
Let $h_{[\BR]}(A)$ be the holonomy of the connection along
the path. This homomorphism of the fiber $P_{\{y, + \infty\}}$ and
$P_{\{y, -\infty\}}$
depends on the equivalence class of $A$, so the construction
defines a map
\[ h_{[\BR]} : \tilde{\CB}_{Y \x {\BR}}^*(a, b) 
\to \mbox{Hom}(P_{\{y, -
\infty\}}, P_{\{y, + \infty\}}) \equiv SO(3) . \]
This map is independent of the time-translation since $h_{[\BR]}(A)
= h_{[\BR]}(A_0)$ for the instanton $A_0 \in \hat{{\CB}}_{Y \x {\BR}} 
(a,b)$. Thus the map factor through $\tilde{\CB}^*_Y$ ( identified
with $\hat{\tilde{\CB}}_{Y \x {\BR}}^*(a,b)$), so  
we obtain $\underline{h_{[\BR]}} : \tilde{\CB}^*_Y \to SO(3)$.
Then we pull-back the fundamental class
$[SO(3)] \in H^3(SO(3),{\BZ})$ to get a cohomology class
$\underline{h_{[\BR]}}^*([SO(3)]) \in
H^3(\tilde{\CB}^*_Y, {\bf Z})$.
By the slant product in \cite{dk} \S 5.1.2, $\mu(y_0 \x {\BR}) = 
\mu (y_0, t_0) = \nu$. See \cite{bd} \S 4.2 for different, but equivalent, 
constructions of this class.

Next we consider the map $\tilde{\mu}_{- \frac{1}{4}p_1}: H_3(Y \x {\BR}) \to 
H^1({\CB}_{Y \x {\CR}}(a, b))$. For the unique fundamental class $[Y]$ of $Y$,
each connection $A$ over $Y \x {\CR}$ can be used to calculate the Chern-Simons
invariant of $A|_{Y \x \{t\}}$. This invariant takes $S^1$-value and 
is independent of the $t$-variable.
Therefore it defines a map $cs_Y: \tilde{\CB}_{Y \x {\CR}}
(a, b) \to S^1$. Hence $cs_Y^*([S^1])$ is an element of 
$H^1(\tilde{\CB}_{Y \x {\CR}}(a, b))$. By the construction
$cs_Y^*([S^1]) \in H^1(\tilde{\CB}_{Y \x {\CR}}(a, b)/{\CR}) = 
H^1(\tilde{\CB}_Y)$. From the Chern-Weil theory of the second Chern class, 
$cs_Y^*([S^1])$ coincides with 
$\tilde{\mu}_{- \frac{1}{4}p_1}([Y]) = \mu(Y)$. So we
have the following.

\begin{pro} \label{2e}
For an integral homology 3-sphere $Y$, 
\[H^*({\CB}_Y; {\BQ}) = H^*(\tilde{{\CB}}_Y; {\BQ}) = 
\text{Sym}^*(\nu) \ox \Lam^*(\mu (Y)).\]
\end{pro}
Proof: By Proposition~\ref{1e} and the discussion above, we have
\[H^*(\tilde{\CB}_Y; {\BQ}) = \text{Sym}^*(\nu) \ox \Lam^*(\mu (Y)).\]
The total space $\tilde{\CB}_Y^*$ of the base-point fibration has the
same weak homotopy type as $\tilde{\CB}_Y$. By the Gysin sequence, we obtain
\[H^*({\CB}_Y^*; {\BQ}) = \oplus \nu^i \cup H^{*-i}(\tilde{\CB}_Y; {\BQ}).\]
So ${\CB}_Y^*$ has the same cohomology ring of $\tilde{\CB}_Y$.
The structure of ${\CB}_Y$ normal to the singular
strata space ${\CB}_Y \setminus {\CB}_Y^*$ 
is a cone on an infinite dimensional
space. The result follows. \qed

Since both classes $\nu$ and $\mu (Y)$ are arisen from the $\mu$-map
construction, we also obtain their Poincar\'{e} dual $V_{y_0}=P.D(\nu)$ and
$V_Y = P.D(\mu (Y))$ in $\tilde{\CB}_Y(a, b)$. The divisors
$V_{y_0}$ and $V_Y$ are $t$-independent so that they also represent
two divisors in the space ${\CB}_Y$.

\subsection{The well-defined action of $H^*({\CB}_Y, {\BQ})$ on $I^{(r)}_*(Y, {\BZ})$}

Let us first recall the Floer-Uhlenbeck compactness on $Y \x {\BR}$.

\begin{df} An {\it ideal anti-self-dual connection (trajectory)} over 
$Y \times {\BR}$, 
of Chern number $k$, is a pair
\[ ( A ; ( x_1,...,x_l ) ) \in {\CM}^{k-l}_{Y \times {\BR}}(a,b) \times
S^l(Y \times {\BR}) \]
where $A$ is a point in ${\CM}^{k-l}_{Y \times {\BR}}(a,b)\cap {\CB}_{\delta }$
and $( x_1,...,x_l )$ is an unordered $l$-tuple of points of $Y \times {\BR}$
\end{df}
Let $\{ A_n\}, n \in {\bf N}$, be a sequence of connections of charge $k$ on the
 $SU(2)$
bundle $P$ over $Y\times {\BR}$. We say that the gauge equivalence classes
$\{A_n \}$ {\it converge weakly} to
 a limiting ideal anti-self-dual connection $(A ; (
x_1,...,x_l ) )$ if
\begin{description}
\item[(i)] The action densities converge as measures, i.e. for any continuous
function $f$ on $Y \times {\BR}$,
\[ \int_{Y \times {\BR}} f |F(A_n)|^2 d \mu \rightarrow \int_{Y \times {\BR}}
 f|F(A)|^2 d \mu  + 8 \pi^2 \sum_{i=1}^l f(x_i).\]
 
\item[(ii)] There are bundle maps  \[\rho_n : P|_{Y\times {\BR} \setminus
\{x_1,...,x_l\} } \rightarrow P|_{Y\times {\BR} \setminus
\{x_1,...,x_l\} } \]
\noindent such that ${\rho }_n^{\ast }(A_n)$ converges to 
$A$ in $C^{\infty }$ on
compact subsets of the punctured manifold.\end{description}
\vspace{0.1 in}
 
\begin{df} Let $a$ and $b$ be flat $SU(2)$ connections over $Y$.
A {\it chain of connections}
$(B_1$, ... , $B_n )$ from $a$ to $b$ is a finite set
of connections over $Y\times {\BR}$
which converge to flat connections $c_{i-1},
c_i$ as $t \to \mp\infty$ such that
$a = c_0$, $c_{n} = b$,
and $B_i$ connects $c_{i-1}, c_i$ for $0 \leq i \leq n$.
\end{df}
We say that the sequence $\{A_{\a }\} \in {\CM}_{Y\times 
{\BR}}^k(a,b)$ is {\it (weakly) convergent} to the chain of connections
$(B_1, ... , B_n )$ if there is a sequence of n-tuples of real numbers 
$\{t_{\a , 1} \leq \cdots \leq t_{\a,n}\}_{\a }$,
such that $t_{\a,i}-t_{\a,i-1}\to\infty$ as $\a \to \infty$, and
if, for each $i$, the translates
$t_{\a,i}^*A_\a=A_{\a }(\circ-t_{\a,i})$ converge weakly to $B_i$.
 
\begin{df} An ``ideal chain connection'' joining flat connections $a$ and $b$
over $Y$ is a set
 \[ (A_j; x_{j 1}, ... ,x_{j l_j} )_{1 \leq j \leq J} \]
where $ (A_j )_{1 \leq j \leq J}$ is a chain connection and for each $j$,
$(A_j; x_{j1}, ... ,x_{jl_j} )$ is an ideal instanton.
\end{df}
In this setup, a version of the Uhlenbeck compactness theorem holds. We state it
in a form
proved by Floer in \cite{fl} (see also in \cite{li1}).
\vspace{0.1 in}
 
\begin{thm}\label{36}[Floer-Uhlenbeck compactness on $Y\times {\BR}$]
Let $A_{\a } \in {\CM}^k_{Y\times {\BR}}\cap {\CB}_{\delta
}(a_{\a
},b_{\a })$  be a sequence of
anti-self-dual connections with uniformly bounded action.
Then there exists a subsequence converging to an ideal chain connection
$(A_j; x_{j 1}, ... ,x_{j l_j} )_{1 \leq j \leq J}$. 
Moreover, one has
\[ \sum_{j=1}^J (k_j + l_j) = k, \hspace{0.3 in}  c_2(A_j) = k_j\;(\mbox{not
necessarily an integer}). \]
\end{thm}

In order to define a cohomology class evaluation on moduli spaces
over $Y \x {\BR}$, we need to rule out the ideal chain connections
over $Y \x {\BR}$. 

\begin{df} \label{action}
For a cohomology class $\o \in H^p({{\CB}}_Y; {\BQ})$, 
the action of $\o$ is given by
\[ \begin{array}{ll}
\o \cap : & C_n^{(r)}(Y) \to C_{n-p-1}^{(r)}(Y) \\
 & a \mapsto \sum_{SF(a^{(r)}) - SF(b^{(r)}) =p+1} \# 
({\hM}_{Y \x {\BR}}(a, b) \cap P.D(\o)) \cdot b,
\end{array} \]
where $\# ({\hM}_{Y \x {\BR}}(a, b) \cap P.D(\o))$ 
is the algebraic number
of $\o$ evaluating on the $p$-dimensional moduli space
${\hM}_{Y \x {\BR}}(a, b)$ in ${{\CB}}_Y$.
\end{df}

Note that every rational cohomology class in $H^*({{\CB}}_Y; {\BQ})$
can be represented as $\o = \nu^k \cup \mu(Y)^l$. So the action of $\o$
defined above is equivalent to the following:
\[(\o \cap )(a) \]\[= \left\{\begin{array}{ll}
\# ({\hM}_{Y \x {\BR}}(a, b) \cap V_{y_1} \cap \cdots V_{y_k} \cap
V_Y^1 \cap \cdots V_Y^l) \cdot b 
& \text{if $SF(a^{(r)}) - SF(b^{(r)}) - 1 = 3k +l$}\\
0 & \text{otherwise.}
\end{array} \right. \]

Since $\nu$ and $\mu (Y)$ are induced from cohomology classes on $4$-manifold
$Y \x {\BR}$, we can use the same construction in \cite{dk} to show that
the intersection
${\CM}_{Y \x {\BR}}(a, b)$ with $V_y$ and $V_Y$ is transversal.
Using the translation invariant, we pass the transversality for the
intersection in ${\CB}_Y$.
By Theorem~\ref{36}, we see that ${\hM}_{Y \x {\BR}}(a, b)$ 
in general does
not represent a genuine cycle in ${{\CB}}_Y$, so the pairing between
${\hM}_{Y \x {\BR}}(a, b)$ 
with cohomology classes in Definition~\ref{action}
needs to be justified.

\begin{lm} \label{nobub}
For $a \in C_n^{(r)}(Y)$ and $b \in C_{n-p-1}^{(r)}(Y)$, the compactification
of ${\hM}_{Y \x {\BR}}(a^{(r)}, b^{(r)})$ 
does not contain any ideal anti-self-dual
connection over $Y \x {\BR}$.
\end{lm}
Proof: Suppose not. There is a sequence 
$\{A_n \in {\CM}_{Y \x {\BR}}(a^{(r)}, b^{(r)})\}$
which converges to an ideal anti-self-dual connection
$(A; (x_1, \cdots, x_i)) \in 
{\CM}^{p+1 - 8i}_{Y \x {\BR}}(\tilde{a}, \tilde{b}) \x S^i(Y \x {\BR})$.
By the definition, we have 
\[cs(a^{(r)}) - cs(b^{(r)}) = \frac{1}{8 \pi^2} \int_{Y \x {\BR}}Tr
F_{A_n} \wedge F_{A_n},  \]
\[ind_{A_n} = SF (a^{(r)}) - SF (b^{(r)}) = p+1.\]
The Chern-Simons is continuous in the sense of weak convergence.
So we get 
\[\lim_{n \to \infty} \frac{1}{8 \pi^2} 
\int_{Y \x {\BR}}Tr F_{A_n} \wedge F_{A_n} =
\frac{1}{8 \pi^2} \int_{Y \x {\BR}}Tr F_{A} \wedge F_{A} + \sum_{j=1}^in_j,\]
where $n_j$ is the multiplicity of $x_j$, and $\sum_{j=1}^in_j \neq 0$ 
by the hypothesis.
For the asymptotic values $\tilde{a}$ and $\tilde{b}$ of $A$, 
we have
\[cs(\tilde{a}) - cs(\tilde{b}) + \sum_{j=1}^in_j = cs(a^{(r)}) - cs(b^{(r)}).\]
If $\sum_{j=1}^in_j \neq 0$, then $\tilde{a} \neq a^{(r)}$ or
$\tilde{b} \neq b^{(r)}$ since $ind_{A} = p+1 - \sum_{j=1}^in_j$.
Thus for the fixed $a^{(r)}$ and $b^{(r)}$, the compactification of 
${\hM}_{Y \x {\BR}}(a^{(r)}, b^{(r)})$ does not contain any ideal anti-self-dual
connection. \qed

Note that the Chern-Simons value is proportional to the spectral flow
by Atiyah-Patodi-Singer's formula. Any ideal anti-self-dual connection
at a point in $Y \x {\BR}$
requires at least energy of spectral flow $8$, in which the Chern-Simons
value decreases its value by $1$. This changes the elements for the different
regular values. So our ${\BZ}$-graded
Floer chain complex rules out this situation.

\begin{lm} \label{comp}
For $\o \in H^p({{\CB}}_Y; {\BQ})$, 
$a \in C_n^{(r)}$ and $b \in C_{n-p-1}^{(r)}$, the intersection
${\hM}_{Y \x {\BR}}(a, b)
\cap P.D(\o)$ is compact.
\end{lm}
Proof: Note that $\o = \nu^k \cup \mu(Y)^l$, so it suffices to show that
the intersection 
${\hM}_{Y \x {\BR}}(a, b) \cap V_{y_1} \cap \cdots V_{y_k}
\cap V_Y^1  \cap \cdots V_Y^l$ is compact.

Suppose $\{A_n\}$ is a sequence in the intersection 
${\hM}_{Y \x {\BR}}(a, b) \cap V_{y_1} \cdots V_Y^l$.
There exists a unique lift of the sequence $\tilde{A}_n \in
{\hM}_{Y \x {\BR}}(a^{(r)}, b^{(r)})$.
By Theorem~\ref{36}, there exists a subsequence converging to an
ideal chain connection. By Lemma~\ref{nobub}, we know that the ideal
anti-self-dual connection does not occur. So
the sequence $A_n$ does not have ideal anti-self-dual connection.
For the ideal
chain connections $(B_1, \cdots, B_n)$ such that
\[B_i \in {\hM}_{Y \x {\BR}}(c_{i-1}^{(r)}, c_i^{(r)}) 
\cap V_{y_1} \cap \cdots V_{y_{k_i}}
\cap V_Y^1  \cap \cdots V_Y^{l_i}, \]
\[\sum_{i=1}^nk_i = k, \ \ \ \ 
\sum_{i=1}^n l_i = l .\]
If $SF (c_{i-1}^{(r)}) - SF (c_i^{(r)}) -1 < 3k_i+l_i$,
by the dimension counting and the general position of the divisors, then 
\[{\hM}_{Y \x {\BR}}(c_{i-1}^{(r)}, c_i^{(r)}) \cap V_{y_1} \cap \cdots 
V_{y_{k_i}}
\cap V_Y^1  \cap \cdots V_Y^{l_i} = \emptyset .\]
So we may have the situation
\begin{equation} \label{big}
SF (c_{i-1}^{(r)}) - SF (c_i^{(r)}) -1 \geq 3k_i+l_i.
\end{equation}
Also we know that
\begin{equation} \label{add}
\sum_{i=1}^n (SF (c_{i-1}^{(r)}) - SF (c_i^{(r)})) = SF(a^{(r)}) -
SF(b^{(r)}) = 3k +l +1.
\end{equation}
Adding (\ref{big}) for $i$, we have
\begin{equation*} \begin{split}
\sum_{i=1}^n (SF (c_{i-1}^{(r)}) - SF (c_i^{(r)}) -1) &= 
\sum_{i=1}^n (SF (c_{i-1}^{(r)}) - SF (c_i^{(r)})) - n \\
& \geq 3 \sum_{i=1}^n k_i + \sum_{i=1}^n l_i \\
& = 3 k + l\\ 
\end{split} 
\end{equation*}
Combining with (\ref{add}), we obtain $n \leq 1$.
So there is no ideal chain connection. The result follows. \qed

\begin{pro} \label{well}
The action $\o \cap$ is well-defined for $\o \in H^*({{\CB}}_Y; {\BQ})$.
\end{pro}
Proof: By Lemma~\ref{comp}, we have 
${\hM}_{Y \x {\BR}}(a, b) \cap V_{y_1} \cap \cdots V_{y_k}
\cap V_Y^1  \cap \cdots V_Y^l$
is a compact, 0-dimensional, oriented manifold. So the algebraic number
in the expression $(\o \cap )(a)$ is finite. It is easy to check that
$\o \cap$ is independent of the choices of the representatives
$V_{y_i} (i = 1, \cdots, k)$ and $V_Y^j (j = 1, \cdots, l)$
(same as in \cite{dk}). Therefore the action $\o \cap$ is well-defined.
\qed

\begin{lm} \label{comu}
The action $\o \cap$ commutes with $\bd^{(r)}$, i.e.,
$\bd^{(r)}_{n-p-1} \circ (\o \cap ) = (\o \cap ) \circ \bd^{(r)}_{n+1}$.
\end{lm}
Proof: For $a \in C_{n+1}^{(r)}(Y)$ and $c \in C_{n-p-1}^{(r)}(Y)$, the space
\[ K = {\hM}_{Y \x {\BR}}(a, c)\cap V_{y_1} \cap \cdots V_{y_k}
\cap V_Y^1  \cap \cdots V_Y^l \]
is a 1-dimensional manifold in ${{\CB}}_Y$. 
The boundary components of $K$ gives
\[({\hM}_{Y \x {\BR}}(a, b) \cap P.D(\o)) \x 
{\hM}_{Y \x {\BR}}(b,c) \coprod
{\hM}_{Y \x {\BR}}(a, b_1) \x 
({\hM}_{Y \x {\BR}}(b_1, c) \cap P.D(\o)).\]
Again the ideal anti-self-dual connections are ruled out by the same method in 
Lemma~\ref{nobub}.
By Lemma~\ref{comp}, the counting argument gives
\[SF(c_{i-1}^{(r)}) - SF (c_i^{(r)}) - 1 \geq 3k_i + l_i, \ \ \ 
3\sum_{i=1}^n k_i + \sum_{i=1}^nl_i = p, \]
\[ \sum_{i=1}^n (SF(c_{i-1}^{(r)}) - SF (c_i^{(r)})) = p+2, \] 
for the ideal chain connection $(B_i)_{1 \leq i \leq n}$. Thus adding the
above inequality we obtain $n \leq 2$. So there is a possible
chain connection $(B_1, B_2) \in {\hM}_{Y \x {\BR}}(a, b) \x
{\hM}_{Y \x {\BR}}(b, c)$ boundary component with
$SF(a^{(r)}) - SF(b^{(r)}) = i$ and $SF(b^{(r)})-SF(c^{(r)}) = p+2 -i$,
where $1 < i < p+1$. Let $d: {\CB}_Y \to {\CB}_Y \x {\CB}_Y$ be the diagonal 
map. So $d^*(\o) \in H^*({\CB}_Y \x {\CB}_Y; {\BQ})$ and
$P.D(d^*(\o)) = P.D(\o) \x P.D(\o) \in H_*({\CB}_Y \x {\CB}_Y; {\BQ})$. 
Thus we have
\[\# ({\hM}_{Y \x {\BR}}(a, b) \x {\hM}_{Y \x {\BR}}(b, c))\cap P.D(d^*(\o)) \]
\begin{eqnarray*}
& = &\# ({\hM}_{Y \x {\BR}}(a, b) \x {\hM}_{Y \x {\BR}}(b, c))\cap
(P.D(\o) \x P.D(\o))\\
& = &\# ({\hM}_{Y \x {\BR}}(a, b)\cap P.D(\o))\x \#
({\hM}_{Y \x {\BR}}(b, c) \cap P.D(\o))\\
& = & 0,
\end{eqnarray*}
by the dimension reason from transversality. Even though there are
boundary components other than we expect, they contribute zero.
Hence the result follows. \qed

The map $\o \cap $ in Definition~\ref{action} defines a 
${\BZ}$-graded chain map, and induces a map (still denoted by
$\o \cap$) on the ${\BZ}$-graded Floer homology:
\[\o \cap : I_n^{(r)}(Y, {\BZ}) \to I_{n-\deg (\o)-1}^{(r)}(Y, {\BZ}).\]
Now for any $\o \in H^*({\CB}_Y; {\BQ})$, there is a well-defined
action on the ${\BZ}$-graded Floer homology. Thus we obtain the desired
cap-product structure
\[H^*({\CB}_Y; {\BQ}) \ox I_*^{(r)}(Y, {\BZ}) \to I_*^{(r)}(Y, {\BZ})\]
on the ${\BZ}$-graded Floer homology.

\noindent{\bf Remark}: It would be nice to show that the cap-product
structure induces a $H^*({\CB}_Y; {\BQ})$-module structure on
$I_*^{(r)}(Y, {\BZ})$. Clearly, we have
$(\o_1 + \o_2) \cap (a) = (\o_1 \cap) (a) + (\o_2 \cap)(a)$, $(1 \cap )(a) = a$
and $(\o \cap)(a_1 + a_2) = (\o \cap)(a_1) + (\o \cap)(a_2)$.
The only thing left is to check the ring homomorphism from
$H^*({\CB}_Y; {\BQ})$ to $\text{End} (I_*^{(r)}(Y, {\BZ}))$.
This is equivalent to $(\o_1 \cup \o_2) \cap = (\o_1 \cap) \circ (\o_2 \cap)$.
Note that $(\o_1 \cup \o_2) \cap$ acts on
$I_n^{(r)}(Y, {\BZ}) \to I_{n-\deg (\o_1) - \deg (\o_2) -1}^{(r)}(Y, {\BZ})$,
but the composition $(\o_1 \cap) \circ (\o_2 \cap)$ acts on
$I_n^{(r)}(Y, {\BZ}) \to I_{n-\deg (\o_2) -1}^{(r)}(Y, {\BZ}) \to
I_{n-\deg (\o_1) - \deg (\o_2) -2}^{(r)}(Y, {\BZ})$.
It is impossible to have 
$(\o_1 \cup \o_2) \cap = (\o_1 \cap) \circ (\o_2 \cap)$. 
The ring structure on $H^*({\CB}_Y; {\BQ})$ is not
preserved by the cap-product in Proposition~\ref{well}. The $-1$ shift is due to
the noncompact ${\BR}$-action on the moduli spaces.

\section{Cap-product structures on the spectral sequence $E^k_{n,j}(Y)$}

In this section, we are going to extend the cap-product structure
to the Fintushel-Stern spectral sequence. The cap-product structure
on $I_*^{(r)}(Y, {\BZ})$ serves as an initial step as the cap-product 
on $E^1_{*, *}(Y)$ in Theorem~\ref{E1}. Then we show that the 
cap-product structure can be deduced on $E^k_{*, *}(Y)$, in particular 
on the ${\BZ}_8$-graded instanton Floer homology.

\begin{lm} \label{41}
For $\o \in H^p({\CB}_Y; {\BQ})$, the action
$\o \cap$ in Definition~\ref{action} is a filtration preserving homomorphism.
\end{lm}
Proof: Recall that the filtration defined in the Fintushel-Stern spectral
sequence is given by
\[F_n^{(r)}C_j(Y) = \sum_{k \geq 0}C_{n+8k}^{(r)}(Y).\]
Thus the map in Lemma~\ref{well} gives a well-defined map
\[\o \cap : C_{n+8k}^{(r)}(Y) \to C_{n+8k-\deg (\o) -1}^{(r)}(Y),\]
which induces a map $\o \cap: F_n^{(r)}C_j(Y) \to F_{n-\deg (\o) -1}^{(r)}
C_{j_{\o}}(Y)$. Here $j_{\o} \equiv n - \deg (\o) - 1 \pmod 8$.
Thus we obtain a filtration preserving map $\o \cap$ with degree 
$ - \deg (\o) - 1$. \qed

Note that $\o \cap$ does not send $C_{n+8k}^{(r)}(Y)$ to 
$C_{n+8(k+1)-\deg (\o) -1}^{(r)}(Y)$. By Lemma~\ref{comu} and Theorem~\ref{E1},
$d^0 = \bd^{(r)}$ in the spectral sequence. So 
\[\o \cap : E^0_{n,j} \to E^0_{n-\deg (\o) -1, j-\deg (\o) -1}(Y) \]
commutes with $d^0$ by Lemma~\ref{comu}. 
There is an induced action (still denoted by $\o \cap$)
\[\o \cap : E^1_{n,j} \to E^1_{n-\deg (\o) -1, j-\deg (\o) -1}(Y).\]
Our spectral sequence is periodical in $j$-direction, so the same map
appears infinitely many times at $j \pmod 8$.
So the rational cohomology ring $H^*({\CB}_Y; {\BQ})$ acts diagonally
on $E^1_{*, *}(Y)$.
The task is to show that this action commutes with all higher differentials
$d^k (k \geq 1)$.
\begin{lm} \label{cmu}
The map $\o \cap$ commutes with $d^k (k \geq 1)$.
\end{lm}
Proof: Note that each element of $E^k_{n, j}(Y)$ is a survivor from previous
differentials, and we only consider the rational cohomology class $\o$.
So $a \in E^k_{n, j}(Y)$ is an element in $F_n^{(r)}C_j(Y)$ with 
$\bd a \in F_{n-1+8k}^{(r)}C_{j-1}(Y)$. The higher differential $d^k$
is again the algebraic counting of 1-dimensional moduli space from
$a \in E^k_{n, j}(Y)$ to $b \in E^k_{n-1 +8k, j-1}(Y)$. 
The following diagram
\begin{equation*}
\begin{CD}
E^k_{n, j} @>d^k>>E^k_{n+8k-1, j-1}\\
 @VV{\o \cap }V @ VV{\o \cap }V\\
E^k_{n-\deg (\o) -1, j-\deg (\o) -1} @>d^k>>E^k_{n-\deg (\o)+8k -2, j-\deg (\o) -2}
\end{CD}
\end{equation*}
is commutative by the same method of proof in Lemma~\ref{comu}. \qed

\begin{thm} \label{main}
There is a well-defined cap-product action of $H^*({\CB}_Y; {\BQ})$ on the
Fintushel-Stern spectral sequence $(E^k_{n,j}(Y), d^k)$.
\end{thm}
Proof: There is a well-defined induced action
$H^*({\CB}_Y; {\BQ}) \ox E^1_{n,j}(Y) \to E^1_{n- \deg (\o) - 1, j- \deg (\o) 
-1}(Y)$ by Lemma~\ref{41}. Inductively, we obtain the cap-production $\o \cap$
on $E^k_{n,j}(Y)$ by Lemma~\ref{cmu} for all $k \geq 1$.
\qed

\begin{cor}
The ${\BZ}_8$-graded instanton Floer homology admits a
$H^*({\CB}_Y; {\BQ})$-group action from the
cap-product structure.
\end{cor}
Proof: By Theorem~\ref{E1}, any one period of the Fintushel-Stern
spectral sequence converges to the ${\BZ}_8$-graded instanton Floer homology.
So the result follows from Theorem~\ref{main} for the
term $E^{\infty}_{n, j}(Y)$.
\qed

\noindent{\bf Remark:}
Note that $H^*({\CB}_Y; {\BQ})$ acts on $E^k_{*, *}(Y)$ as a group, 
not ring. So there is no 
$H^*({\CB}_Y; {\BQ})$-module structure from our definition of
the action.

\noindent{\bf Example:} Let $Y = \Si (2, 3, 5)$ be the Poincar\'{e} 3-sphere.
By a result of Fintushel-Stern \cite{fs}, 
\[I^{(0)}_1(Y, {\BZ}) = {\BZ}\la a_{\a} \ra;
I^{(0)}_5(Y, {\BZ}) = {\BZ}\la a_{\b} \ra;
I^{(0)}_i(Y, {\BZ}) = 0, \ \ i \neq 1, 5.\]
By a result of Kronheimer (see also \cite{au} Lemma 5.1), we have
\[\nu \cap : C_5^{(0)}(Y, {\BZ}) \to C_1^{(0)}(Y, {\BZ}) \]
given by $(\nu \cap )(a_{\b}) = 2 a_{\a}$, and $\mu(Y) \cap $ acts trivially 
by the degree reason. Any other cohomology class $\nu^k \cup \mu(Y)^l (l \neq 0)$ acts trivially on $I_*^{(0)}(Y, {\BZ})$. In this case, the Fintushel Stern
spectral sequence collapses at $E^1_{*, *}(Y)$. So the
cap-product structure on the ${\BZ}_8$-graded instanton
Floer homology has only one nontrivial action
\[\nu \cap : HF_5(Y, {\BZ}) \to HF_1(Y, {\BZ}).\]

\end{document}